\documentclass[graybox]{svmult}
\usepackage{amsmath}
\usepackage{amsfonts}
\usepackage{amssymb}
\newcommand{\Section}[1]{\section{#1} \setcounter{equation}{0}}

\DeclareMathOperator{\card}{card}

\makeatletter
\def\imod#1{\allowbreak\mkern6mu({\operator@font mod}\,\,#1)}
\makeatother

\begin{document}

\title{The generalized Montgomery-Hooley formula: A survey}
\author{R. C. Vaughan }
\titlerunning{Generalized Montgomery-Hooley formula}
\authorrunning{R. C. Vaughan}
\institute{R. C. Vaughan \at Department of Mathematics, Penn. State University, University Park, PA 16802, USA\\ \email{rcv4@psu.edu} }

\maketitle

\begin{abstract}
This memoir is a survey of theorems and inequalities which have grown out of, and extended, the seminal estimate of Montgomery \cite{HM70}
\begin{multline*}
V(x,Q)=\sum_{q\le Q}\sum_{\substack{a=1\\
(a,q)=1}}^q \left|
\psi(x;q,a) - \frac{x}{\phi(q)}
\right|^2 \\
= Qx\log x + \textstyle O\big(Qx\log\frac{2x}{Q}\big) + O\big(x^2(\log x)^{-A}\big).,
\end{multline*}
\end{abstract}

\Section{Origins}
\label{sec:one}

\noindent This survey is of necessity eclectic.  It is an expanded version of a talk given in the Number Theory Web Seminar in May 2022 and largely reflects those aspects of the topic which have particularly interested me.  I apologise to anyone whose work has been overlooked.
\par
The eponymous theorem which is the starting point of this article is the following.
\begin{theorem}[Montgomery,\cite{HM70}]
\label{thm:one1}
Define
\begin{equation}
\label{eq:one1}
V(x,Q)=\sum_{q\le Q}\sum_{\substack{a=1\\
(a,q)=1}}^q \left|
\psi(x;q,a) - \frac{x}{\phi(q)}
\right|^2,
\end{equation}
where as usual
\[
\psi(x;q,a)=\sum_{\substack{n\le x\\
n\equiv a\imod q}} \Lambda(n).
\]
Let $A>0$ and suppose that $x>x_0(A)$.  Then for $Q\le x$,
\[
V(x,Q) = Qx\log x + \textstyle O\big(Qx\log\frac{2x}{Q}\big) + O\big(x^2(\log x)^{-A}\big).
\]
\end{theorem}
As with all theorems of this kind there are concomitant conclusions when the function $\psi$ is replaced by $\vartheta$ or $\pi$.  For consistency I will persist with $\psi$.
\par
Montgomery's theorem was refined, with a much simpler proof, in
\begin{theorem}[Hooley,\cite{CH75a}]
\label{thm:one2}
On the hypothesis of the previous theorem,
\[
V(x,Q) = Qx\log Q -cQX + O\big(Q^{\frac54}x^{\frac34} + x^2(\log x)^{-A}\big).
\]
\end{theorem}
One immediate observation.  The error terms of these theorems become less precise when $Q$ is close to $x$.  There is a good reason for this.  When $q\approx x$ the number of residue classes is greater than the number of primes, so $\frac{x}{\phi(q)}$ has to be a bad approximation.
\par
Earlier, Barban \cite{MB64} (see also Barban \cite{MB63}, \cite{MB66}) had established that if $B=B(A)$, $x>x_0(A)$, $Q\le x(\log x)^{-B}$,  then
\[
V(x,Q)\ll x^2(\log x)^{-A}
\]
and this was refined by Davenport and Halberstam \cite{DH66} with $B=A+5$ and Gallagher \cite{PG67} with $B=A+1$.  Also Barban had apparently stated that
\[
V(x,x) = x^2\log x -cx^2 + O\big(x^2(\log x)^{-A}\big).
\]
I have not seen this paper and it would be interesting to see what proof Barban had in mind.
\par
These results perhaps are not very surprising.  After all if one averages over enough things one should be able to establish a precise conclusion.  On the other hand they say that on average
\[
\psi(x;q,a)-\frac{x}{\phi(q)} \ll x^{1/2}q^{-1/2}(\log x)^{1/2}
\]
which is stronger than the generalised Riemann hypothesis.
\par
Thus it is not without some interest to try to understand what ingredients are necessary for success and the extent to which they can be applied.  Gallagher's proof is instructive because it reveals some of those ingredients.
\par
This begins by using Dirichlet characters to pick out the residue classes, and then applies orthogonality followed by the prime number theorem to deal with $\chi_0$.  This leads to
\[
\sum_{q\le Q} \frac1{\phi(q)}\sum_{\substack{\chi\imod q\\
\chi\not=\chi_0}} |\psi(x;\chi)|^2.
\]
Then in this one replaces each character $\chi$ by the primitive character $\chi^*$ of conductor $r$ which induces $\chi$.  This leads to essentially
\[
\sum_{m\le Q}\frac1{\phi(m)}\sum_{1<r\le Q/m}\frac1{\phi(r)} \sum_{\chi^*\imod r} |\psi(x;\chi^*)|^2.
\]
Partial summation and the large sieve then gives
\[
\sum_{L<r\le Q/m}\frac1{\phi(r)} \sum_{\chi^*\imod r} |\psi(x;\chi^*)|^2 \ll \left(
\frac{x}{L}+\frac{Q}{m}
\right)x\log x
\]
where $L=(\log x)^B$, say.  The final ingredient is the Siegel-Walfisz theorem to cover the $r\le L$.
\par
How about the asymptotic formula?  An obvious line of attack is to square out and add in any non-zero $\psi(x;q,a)$ with $(q,a)>1$.  Then
\[
V(x,Q) = S_0+2S_1-2S_2+S_3 +O\big((Q+x)\log^4 x\big)
\]
where
\[
S_0=Q\sum_{n\le x}\Lambda(n)^2,
\]
\[
S_1=\sum_{q\le Q}\sum_{\substack{m<n\le x\\
m\equiv n\imod q}} \Lambda(m)\Lambda(n),
\]
\[
S_2=\sum_{q\le Q} \psi(x)\frac{x}{\phi(q)},\quad S_3=\sum_{q\le Q} \frac{x^2}{\phi(q)}.
\]
All the sums here are easy to deal with except $S_1$.  If one can obtain asymptotic expressions for each sum, then one might expect that the main terms will largely cancel and an approximation for V(x,Q) will drop out from lower order terms.
\par
Montgomery's proof deals with
\[
S_1=\sum_{q\le Q} \sum_{\substack{m<n\le x\\
m\equiv n\imod q}} \Lambda(m)\Lambda(n),
\]
by writing this as
\[
\sum_{h\le x} d_Q(h) R(x;h)
\]
where
\[
d_Q(h) =\sum_{\substack{q|h\\
q\le Q}} 1,\quad R(x;h) = \sum_{\substack{m,n\le x\\
n-m=h}} \Lambda(m)\Lambda(n).
\]
One then appeals to Vinogradov's method in additive prime number theory to replace $R(x;h)$ by $\mathfrak S(h)(x-h)$ where $\mathfrak S$ is the appropriate singular series.  The relevant theorem here actually is due to Lavrik \cite{AL60}.  Although superseded by Hooley's idea I will return to this later.
\par
Hooley's idea is as follows.  Let the large sieve deal with $q\le Q_0=x(\log x)^{-B}$ and suppose $Q_0<Q\le x$.
We can then take the difference of two sums of the kind
\[
S_1'(Q)=\sum_{Q<q\le x}\sum_{\substack{m<n\le x\\
m\equiv n\imod q}} \Lambda(m)\Lambda(n).
\]
Write this as
\begin{align*}
\sum_{m<n\le x} \sum_{\substack{qr=n-m\\
Q<q\le x}}\Lambda(m)\Lambda(n)&= \sum_{m<n\le x} \sum_{\substack{r|n-m\\
r<\frac{n-m}{Q}}}\Lambda(m)\Lambda(n) \\
&= \sum_{r<\frac{x}{Q}} \sum_{m< x-rQ}\Lambda(m) \sum_{\substack{
m + rQ < n\le x\\
r|n-m
}}\Lambda(n).
\end{align*}
Now this Hooley inversion enables the proof to be completed by a simple aplication of Siegel-Walfisz.
\par
In this survey I am not so concerned with refining these results, or speculation about sums like
\[
\sum_{\substack{a=1\\
(a,q)=1}}^q\left|
\psi(x;q,a)-\frac{x}{\phi(q)}
\right|^2
\]
although these are of fundamental interest.  My main concern is the extent to which these ideas can be applied to functions significantly different from $\psi(x;q,a)$.  Hooley wrote at least 20 papers in some of which these ideas are extended to a wide class of functions.  The primes have the advantage that they are uniformly distributed into the reduced residue classes.  Most sequences of number theoretic interest are not so well behaved.  Even the square free numbers are deficient in this regard.

\Section{Hooley III}
\label{sec:two}

\noindent The first paper which looks at a general class of cognate problems is Hooley III \cite{CH75c}. He requires an analogue of the Siegel-Walfisz theorem which, for some unfathomable reason,
is labelled
\par
{\bf Criterion U}.  {\it Let $\mathcal S\subset\mathbb N$ and
\[
S(x;q,a)=\sum_{\substack{s\in\mathcal S, s\le x\\
s\equiv a\imod q}} 1
\]
and suppose that for $x>x_0(A)$ we have
\[
S(x;q,a)=f\big(q,(q,a)\big) x + O\big(x(\log x)^{-A}\big).
\]
}
The dependence of the main term on $(q,a)$ rather than $a$ is satisfactory for many applications, such as the squarefree numbers, but nevertheless signals a dependence on the large sieve.  His final conclusion is that, when $Q\le x$,
\[
\sum_{q\le Q}\sum_{a=1}^q |S(x;q,a)-f\big(q,(q,a)\big) x|^2 \ll Qx + x^2(\log x)^{-A}.
\]
There are a number of other generalizations of these techniques.  For example Smith \cite{ES10} has established a version of the Montgomery-Hooley theorem when $\psi(x;q,a)$ is replaced by
\[
\theta_K(x;q,a) =\sum_{\substack{\mathrm N\mathfrak p\le x\\
\mathrm N\mathfrak p\equiv a\imod q}} \log \mathrm N\mathfrak p
\]
on the assumption that $K$ is a Galois extension of $\mathbb Q$.
\par
Another example is due to Keating and Rudnick \cite{KR14}.  There they establish an analogue of Montgomery-Hooley for function fields.

\Section{The Hardy-Littlewood method}
\label{sec:five}

\noindent At this point it is useful to introduce another perspective on the methods so far discussed.
A significant proportion of the work in the area had been consequent on the assumption of the generalised Riemann Hypothesis, and in Goldston and Vaughan  \cite{GV96} an idea was introduced which, whilst facilitating the use of that hypothesis, might be thought of as being a backwards step.  However it transpires that it plays a signicant r\^ole in some later work, and it permitted the following theorem to be established.
\begin{theorem}[Goldston\& Vaughan  \cite{GV96}]
\label{thm:two1}
Suppose that the generalised Riemann hypothesis holds and let
\begin{equation}
\label{eq:two1}
U(x,Q)=V(x,Q)-Qx\log Q-cxQ
\end{equation}
where $V(x,Q)$ is given by (\ref{eq:one1}), and
\begin{equation}
\label{eq:two2}
c=\gamma+\log(2\pi)+1+\sum_p\frac{\log p}{p(p-1)}.
\end{equation}
Then (i) when $1\le Q\le x$ one has
\begin{equation}
\label{eq:two3}
U(x,Q)\ll Q^2(x/Q)^{1/4+\epsilon}+x^{3/2}(\log 2x)^{5/2}(\log\log
3x)^2.
\end{equation}
and (ii) there is an absolute constant $C$ such that when
$x/Q\rightarrow\infty$ with
$$Cx^{5/7}(\log 2x)^{10/7}(\log\log 3x)^{8/7}<Q\le x$$
one has
\begin{equation}
\label{eq:two4}
U(x,Q)Q^{-2}=\Omega_\pm\left(
(x/Q)^{1/4}\right).
\end{equation}
\end{theorem}
\par
Recall that a key ingredient to the original result is an estimate for
\[
S_1=\sum_{m<n\le x} \sum_{\substack{qr=n-m\\
q\le Q}} \Lambda(m)\Lambda(n)
\]
and that Montgomery's original method was based on Vinogradov's method.  Instead one can write directly
\begin{equation}
\label{eq:two5}
S_1=\int_0^1 F(\alpha)|G(\alpha)|^2 d\alpha,
\end{equation}
where
\begin{equation}
\label{eq:two6}
F(\alpha) = \sum_{q\le Q}\sum_{r\le x/q} e(\alpha qr)
\end{equation}
and
\begin{equation}
\label{eq:two7}
G(\alpha)=\sum_{n\le x} \Lambda(n)e(\alpha n).
\end{equation}
The sum $F$ is essentially trivial to estimate on standard minor arcs, so one can avoid Vinogradovs's method.  Whilst not as simple as Hooley's, it has some advantage of flexibility and avoids the large sieve.  Thus it opens up the possibility of dealing with sequences which are not so well distributed.
\par
This was exploited in Vaughan  \cite{RV98a}, \cite{RV98b}.   What is interesting for us is that this can be pushed further to obtain Montgomery-Hooley style asymptotics for the subject of Hooley III.  It is also clear that the underlying ideas work just as well with $S(x;q,a)$ replaced by
\begin{equation}
\label{eq:two09}
A(x;q,a)=\sum_{\substack{n\le x\\
n\equiv a\pmod q}} a_n
\end{equation}
where the sequence $\{a_n\}$ has the potential to be quite general.  We will suppose it is a real sequence as that captures most interesting examples, but there is no reason in principle why it should not be complex.
\par
We are now concerned with the variance
\begin{equation}
\label{eq:two10}
V(x,Q)=\sum_{q\le Q} \sum_{a=1}^q \left|
A(x;q,a)-f(q,(q,a))x\right|^2
\end{equation}
where $f$ appropriately reflects the local properties of the sequence $\{a_n\}$.
\par
Criterion U can be replaced by a more general assumption, namely that there is an increasing function $\Psi(x)$, with
$\Psi(x)>\log x$ for all large $x$, $\Psi(1)>0$ and
$\int_1^x\Psi(y)^{-1}dy\ll x\Psi(x)^{-1}$, such that
\begin{equation}
\label{eq:two11}
A(x;q,a)=xf(q,(q,a))+ O\left(
\frac{x}{\Psi(x)}\right)
\end{equation}
uniformly for all real $x\ge1$ and natural numbers $q$ and $a$, and we note that immediately from these assumptions we have
$\Psi(x)\ll x$.\par
The most natural assumption concerning the $a_n$ is not that it be the indicator function
of a set but rather that it be bounded in mean square, or, more
precisely, that
\begin{equation}
\label{eq:two12}
\sum_{n\le x}a_n^2\ll x
\end{equation}
uniformly for all positive real $x$.  An important r\^ole is then played by the function
\begin{equation}
\label{eq:two13}
g(q)=\phi(q)\left(
\sum_{r|q}f(q,r)\mu(q/r)\right)^2
\end{equation}
One consequence of (\ref{eq:two12}) is that the series
$$\sum_{q=1}^\infty g(q)$$
converges, and the quality of the main conclusions depends on the rate of convergence of this series and the extent to which
$$x\sum_{q=1}^\infty g(q)$$
is a good approximation to the left hand side of (\ref{eq:two12}).
\par
With the above definitions it is possible to state a simple conclusion.
\bigskip
\begin{theorem}[Vaughan  \cite{RV98b}, Theorem 1]
\label{thm:two2}
Suppose
\begin{equation}
\label{eq:two14}
Q>\sqrt x\log 2x,
\end{equation}
(\ref{eq:two12}) holds, and let
\begin{equation}
\label{eq:two15}
E(z)=\int_1^z\sum_{q>y}g(q)dy
\end{equation}
and
\begin{equation}
\label{eq:two16}
U(x,Q)=V(x,Q)-Q\sum_{n\le x}a_n^2+Qx\sum_{q=1}^\infty g(q).
\end{equation}
Then $E(z)=o(z)$ as $z\rightarrow\infty$ and
\begin{multline*}
U(x,q)\ll x^{3/2}\log x + x^2(\log2x)^{9/2}\Psi(x)^{-1} \\
+ x^2(\log x)^{4/3}\Psi(x)^{-2/3} +
Q^2E\left(\frac{x}{Q}\right).
\end{multline*}
\end{theorem}
\par
By Parseval's identity
$$\int_0^1|G(\alpha)|^2d\alpha=\sum_{n\le x}a_n^2$$
where now
\begin{equation}
\label{eq:two17}
G(\alpha)=\sum_{n\le x}a_ne(n\alpha)
\end{equation}
and it is not hard to show on the assumption (\ref{eq:two11}) that the
contribution from the consequential natural major arcs is asymptotically
$$x\sum_{q=1}^\infty g(q).$$
Hence the main term
$$Q\sum_{n\le x}a_n^2-Qx\sum_{q=1}^\infty g(q)$$
in Theorem \ref{thm:two2} is closely related to the minor arcs.  In many of the common
situations matching our conditions it is known that the contribution
from the minor arcs is smaller than that from the major arcs.  For
example, this is so when $a_n$ is the indicator function of the
$k$-free numbers ($k\ge2$).  Thus, in such a situation the two
expressions in the main terms are largely cancelling.  However, we
can then anticipate that provided we have some knowledge of the
asymptotic behaviour of their difference, and perhaps also of $E(y)$,
it is still possible to obtain the asymptotic behaviour of $V(x,Q)$.
That further information regarding $E(y)$ may be helpful is born out
by the case of $k$-free numbers, for which see \S\ref{sec:three} below, where the final main term is indeed
of the same order of magnitude as $Q^2E(x/Q)$ for a large range of Q.
\begin{theorem}[Vaughan  \cite{RV98b}, Theorem 2]
\label{eq:two17a}
Suppose
\begin{equation}
\label{eq:two18}
\sum_{n\le x}a_n^2-x\sum_{q=1}^\infty g(q) = o\left(
x^{\frac{2+\eta}{2+2\eta}}\right)
\end{equation}
as $x\rightarrow\infty$ and
\begin{equation}
\label{eq:two19}
\sum_{q>y}g(q)\sim cy^{-\eta}
\end{equation}
as $y\rightarrow\infty$ where $\eta$ and $c$ are positive real numbers with $0<\eta<1$.  Suppose further that
$$Q>\sqrt x\log 2x$$
Then
\begin{multline*}
V(x,Q)=Q^2M(x/Q)\\
+O\left(
x^{\frac32}\log x + x^2(\log2x)^{\frac92}\Psi(x)^{-1} + x^2(\log
x)^{\frac43}\Psi(x)^{-\frac23} \right)
\end{multline*}
where
$$M(y)\sim c\frac{-2\zeta(-\eta)}{1-\eta^2}y^{1-\eta}$$
as $y\rightarrow\infty$.
\end{theorem}
The method of proof of Theorem 2 is equally valid under more general
conditions than (\ref{eq:two18}) and (\ref{eq:two19}).  For example, with appropriate
adjustments to (\ref{eq:two18}) and the conclusion, the condition (\ref{eq:two19}) could
be replaced by
\begin{equation}
\label{eq:two20}
\sum_{q>y}g(q)\sim \kappa(y)
\end{equation}
where $\kappa(y)$ is a suitably smooth function with
\[
\lim_{y\rightarrow\infty}\kappa(y)=0.
\]
\par
It is natural to ask whether the main terms in Theorem \ref{thm:two2} are always cancelling and it can be shown that this
is not so by the construction of an example.  The point is that the
example places a positive proportion of the mass in
$$\int_0^1|G(\alpha)|^2d\alpha$$
on the minor arcs.
\begin{theorem}[Vaughan  \cite{RV98b}, Theorem 3]
\label{thm:two4}
Suppose
\[
\lambda=\frac{\sqrt5-1}{2},\text{ and }\theta\in(0,1),
\]
and let $a_n$
be $1$  when $\{\lambda n\}<\theta$ and be $0$ otherwise.  Then
(\ref{eq:two11}) holds with $f(q,(a,q))=\theta/q$ and $\Psi(x)=x^{1/3}$, and
$$\sum_{q=1}^\infty g(q)=\theta^2$$
but
$$\sum_{n\le x} a_n^2=\theta x + O(x^{2/3}).$$
\end{theorem}

\Section{The Squarefree numbers}
\label{sec:three}

\noindent In the context of this survey, the squarefree numbers have a substantial history, Warlimont \cite{RW69}, \cite{RW72}, \cite{RW80}, Orr \cite{RO69}, \cite{RO71}, Croft \cite{MC75}, Vaughan  \cite{RV05}, and Parry \cite{TP21}.
\par
Let $\mu_k$ be the indicator function of the $k$-free numbers,
\begin{equation}
\label{eq:three1}
Q_k(x;q,a)=\sum_{\substack{n\le x\\
n\equiv a\imod q}}\mu_k(n),
\end{equation}
\begin{equation}
\label{eq:three2}
f(q,a) = \sum_{\substack{m=1\\
(m^k,q)|a}}^{\infty} \frac{\mu(m)(m^k,q)}{m^kq},
\end{equation}
\begin{equation}
\label{eq:three3}
V(x,Q)=\sum_{q\le Q} \sum_{a=1}^q |Q_k(x;q,a)-xf(q,a)|^2.
\end{equation}
\begin{theorem}[Vaughan  \cite{RV05}]
\label{thm:three1}
There are $c_k^*>0$, $\tilde c_k>0$ so that, whenever $Q\le x$, we have
\begin{multline}
\label{eq:three11}
V(x,Q) =  c_k x^{\frac1k}Q^{2-\frac1k}\\
+ O\left(
x^{\frac1{2k}}Q^{2-\frac1{2k}}\exp\left(
-c_k^*\frac{(\log 2x/Q)^{\frac35}}{(\log\log 3x/Q)^{\frac15}}
\right)
\right)\\
+ O\left(
x^{1+\frac1k}\exp\left(
-\frac{\tilde c_k(\log x)^{\frac35}}{(\log\log x)^{\frac15}}
\right)
\right)
\end{multline}
where
\begin{equation}
\label{eq:three12}
c_k=\frac{2C_1k^2\big(-\zeta(1/k-1)\big)}{\zeta(k)^2(2k-1)(k-1)},\quad C_1 = \prod_p \left(
1 +
\frac{\nu(p)}
{p^2(p^k-1)^2}
\right)
\end{equation}
and
\[
\nu(p) = (p-1)^2\sum_{j=1}^{k-1}p^{(2-\frac1k)j} - 2p^{2k}+p^{2k-1}+2p^{k+1}-p
\]
\end{theorem}
The error terms here are related to the zero free region of $\zeta$ only, and {\it not} other Dirichlet $l$-functions.
\par
The method used is just as described above, but with great care taken over the consequent main term.  Once more one can note the lack of uniformity as $Q\rightarrow x$ even though now the $k$-free numbers have positive density.

\Section{Hooley VIII}
\label{sec:six}

\noindent Let me now advert to another Hooley paper, Hooley VIII \cite{CH98a}.  Here he deals with a problem which is not directly of the kind which is our central interest, namely a third moment
\[
\sum_{q\le Q} \phi(q) \sum_{\substack{a=1\\
(a,q)=1}}^q \left(
\psi(x;q,a) - \frac{x}{\phi(q)}
\right)^3
\]
Note the weight $\phi(q)$.  This is necessary, since one is expecting that the sum over $a$ is behaving roughly like
\[
x^{3/2}\phi(q)^{-1}
\]
and so for the raw sum we could expect to get no strong beneficial effect from the larger $q$ when we sum over $q\le Q$.  Thus some weight which emphasises the larger $q$ is highly desirable.  The $\phi(q)$ seems rather unnatural compared with a smooth weight, but it is there to alleviate some of the not inconsiderable difficulties that Hooley runs into.
\par
Hooley's method is to follow the pattern established for the second moment.
\[
\sum_{q\le Q} \phi(q) \sum_{\substack{a=1\\
(a,q)=1}}^q \left(
\psi(x;q,a) - \frac{x}{\phi(q)}
\right)^3
\]
Thus the cube is multiplied out and four sums are obtained.  Then asymptotic formulae are established  for each one.
The hardest, coming from the product of three von Mangoldt functions, can be dealt with by Vinogradov's method.
However the really big problem is to show that the main terms sum to $0$.  This is a major achievement and takes many pages.  It also results in the paper being littered with quotations from Dante's Inferno.

\Section{A novel main term}
\label{sec:seven}

\noindent The use of the Hardy-Littlewood method as in Goldston and Vaughan suggests a way of simplifying the difficulties in Hooley VIII.  Recall from (\ref{eq:two5}), (\ref{eq:two6}), (\ref{eq:two7}) that the core problem for the primes concerns
\[
S_1=\int_0^1 F(\alpha)|G(\alpha)|^2 d\alpha.
\]
For $\Lambda$ we expect that on the major arcs, say $\alpha$ with
\[
\textstyle |\alpha -\frac{b}{r}|\le \frac{(\log x)^B}{rx}, \, 1\le b\le r\le R=(\log x)^B, \, (r,b)=1,
\]
\[
G(\alpha)\sim G^*(\alpha) = \sum_{r\le R}\sum_{\substack{b=1\\
(r,b)=1}}^r\frac{\mu(r)}{\phi(r)}\sum_{n\le x} e\big((\alpha-b/r)n\big).
\]
If we rewrite the right hand side as
\[
\sum_{n\le x} e(\alpha n) \sum_{r\le R} \frac{\mu(r)}{\phi(r)} c_r(n)
\]
it makes sense to replace the approximation $x/\phi(q)$ by
\[
\sum_{\substack{n\le x\\
n\equiv a\imod q}} \sum_{r\le R} \frac{\mu(r)}{\phi(r)} c_r(n).
\]
\par
Let
\[
\Xi_R(n) = \sum_{r\le R} \frac{\mu(r)}{\phi(r)} c_r(n),
\]
\[
\rho(x;q,a) = \sum_{\substack{n\le x\\ n\equiv a\imod q}} \Xi_R(n).
\]
Then we can establish the following theorems.
\begin{theorem}[Vaughan \cite{RV03a} Corollary 4.1]
\label{thm:seven1}
If $x>x_0(A)$, $Q\le x$ and $R\le (\log x)^A$, then
\begin{multline*}
\sum_{q\le Q} \sum_{a=1}^q\big(\psi(x;q,a)-\rho(x;q,a)\big)^2 \\
= Qx(\log x/R) -cQx + O\big(
QxR^{-1/2}+x^2(\log x)^2R^{-1}
\big).
\end{multline*}
\end{theorem}
The right hand side has the remarkable feature that when $2\pi x/R<Q\le x$, the main term is smaller than that in Theorems \ref{thm:one1} and \ref{thm:one2} and the error term is uniform in $Q$.
\begin{theorem}[Vaughan \cite{RV03b} Theorem 8]
\label{thm:seven2}
Suppose that
\[
R\le (\log x)^A.
\]
Then
\begin{multline*}
\sum_{q\le Q} q \sum_{a=1}^q \big(\psi(x;q,a)-\rho(x;q,a)\big)^3 \\
=\frac12Q^2x(\log x)^2 + O\big(x^3(\log x)^5R^{-1} +  Q^2x(\log x)\log R \big).
\end{multline*}
\end{theorem}
Note also better than square root cancellation as well as uniformity as $Q\rightarrow x$ in the second theorem.
\par
The simplest proof of
\[
\sum_{q\le Q} \sum_{a=1}^q\big(\psi(x;q,a)-\rho(x;q,a)\big)^2
\]
\[
= Qx(\log x/R) -cQx + O\big(
QxR^{-1/2}+x^2(\log x)^2R^{-1}
\big).
\]
is perhaps still by Hooley's inversion method.  However one can write
\[
\psi(x;q,a)-\rho(x;q,a) = \sum_{\substack{n\le x\\
n\equiv a\imod q}} (\Lambda(n)-\Xi_R(n))
\]
and so
\begin{multline*}
\sum_{q\le Q} \sum_{a=1}^q\big(\psi(x;q,a)-\rho(x;q,a)\big)^2
=  \sum_{n\le x} (\Lambda(n)-\Xi_R(n))^2 \\
+ 2\int_0^1 F(\alpha)|G(\alpha)-G^*(\alpha)|^2 d\alpha.
\end{multline*}
The Hardy-Littlewood method applies directly to the integral and shows that it is small compared with the main term.
\par
The treatment of
\[
\sum_{q\le Q} q\sum_{a=1}^q\big(\psi(x;q,a)-\rho(x;q,a)\big)^3
\]
is to write $\Delta(n) = \Lambda(n) - \Xi_R(n)$
\[
\psi(x;q,a)-\rho(x;q,a) = \sum_{\substack{n\le x\\
n\equiv a\imod q}} \Delta(n)
\]
and cube it out.  Then the core part is
\[
\sum_{q\le Q} \sum_{r,s} \sum_{\substack{l<m<n\le x\\
m-l=qr,\,n-m=qs}} \Delta(l)\Delta(m)\Delta(n)
\]
Let $E(\theta)=G(\theta)-G^*(\theta)$.  Then this can be written as
\[
\int_0^1\int_0^1 F(\alpha,\beta) E(\alpha)E(\beta-\alpha)E(-\beta)d\alpha d\beta
\]
where now
\[
F(\alpha,\beta) = \sum_{q\le Q}\sum_{r\le x/q}\sum_{s\le x/q} e(\alpha qr+\beta qs)
\]
Again the Hardy-Littlewood method is amenable.

\Section{Bad behaviour}
\label{sec:eight}

\noindent In all of the cases so far
\[
\sum_{n\le x} a_n
\]
and the approximations to
\[
\sum_{n\le x, q|n-a} a_n
\]
and
\[
\sum_{n\le x} a_n^2
\]
are well behaved.
\par
Dancs \cite{MD02} has pushed the envelope by considering
\begin{equation}
\label{eq:eight1}
V(x,Q) = \sum_{q\le Q} \sum_{a=1}^q \big(A(x;q,a) - \pi x f(q,a)\big)^2
\end{equation}
with
\begin{equation}
\label{eq:eight2}
a(n)=r(n)=\card\{x,y\in\mathbb Z^2:x^2+y^2=n\}
\end{equation}
and
\begin{align*}
f(q,a) &= q^{-2}\card\{x,y\in\mathbb Z_q:x^2+y^2\equiv a\imod q\} \\
&=q^{-3}\sum_{b=1}^q S(q,b)^2e(-ab/q)
\end{align*}
where
\[
S(q,b)=\sum_{x=1}^q e(bx^2/q).
\]
For $Q\le x$ he obtains
\begin{equation}
\label{eq:eight3}
V(x,Q) = 8Qx\big(\log(x/Q) + C_1\big) + 4Q^2\log Q + C_2 Q^2 + O(x^{5/3+\varepsilon}).
\end{equation}
\par
A trickier example which was looked at by Motohashi \cite{YM73} in the special case $Q=x$ and in the general case by Pongsriiam \cite{PP12} (see also Pongsriiam and Vaughan  \cite{PV18}, \cite{PV21}) is
\begin{equation}
\label{eq:eight4}
A(x;q,a)=\sum_{\substack{n\le x\\
n\equiv a\imod q}} d(n).
\end{equation}
One significant problem here is that in the approximation the local and global factors do not split.  The most useful approximation is
\begin{equation}
\label{eq:eight5}
M(x;q,a) = \frac{x}{q} \sum_{r|q} \frac{c_r(a)}r \left(
\log\frac{x}{r^2} +2\gamma -1
\right).
\end{equation}
Now let
\begin{equation}
\label{eq:eight6}
V(x,Q) = \sum_{q\le Q} \sum_{a=1}^q \big(A(x;q,a) - M(x;q,a)\big)^2.
\end{equation}
When $Q\le x$ he obtains
\begin{multline}
\label{eq:eight7}
V(x,Q) = \frac{Qx}{\pi^2}\left(
\log\frac{Q^2}{x}
\right)^3 +\\
QxP(\log x,\log Q) + O\left(
x^{\frac74}(\log x)^3 + Q^2(x/Q)^{\varepsilon}(\log x)^2
\right)
\end{multline}
where $P(\xi,\eta)$ is a polynomial of degree $2$.   Probably the error terms in (\ref{eq:eight3}) and (\ref{eq:eight7}) are susceptible to some improvements.
\par
Another curious example has been studied by Penyong Ding \cite{PD21}.  He considers
\begin{equation}
\label{eq:eight8}
A(x;q,a)=\sum_{\substack{n\le x\\
n\equiv a\imod q}} r_3(n)
\end{equation}
where
\begin{equation}
\label{eq:eight9}
r_3(n) = \card\{l_1,l_2,l_3\in\mathbb N^3: l_1^3+l_2^3+l_3^3\}
\end{equation}
and uses the approximation
\[
\Gamma(4/3)^3x\rho(q,a)q^{-3}
\]
where $\rho(q,a)$ is the number of solutions of the congruence
\[
l_1^3+l_2^3+l_3^3\equiv a\imod q.
\]
\par
The function $r_3(n)$ is somewhat mysterious since we don't know how
\[
\sum_{n\le x} r_3(n)^2
\]
behaves.  The best that we know (Vaughan \cite{RV21}, Corollary 1.3) is that
\[
x\ll \sum_{n\le x} r_3(n)^2 \ll x^{\frac76}(\log x)^{\varepsilon -\frac52}.
\]
We might expect that it is $\sim cx$, but Hooley \cite{CH86} has shown that if true the value of $c$ is not obvious.
\par
For $r_3(n)$ it is natural to approach the question by using a variant of the Hardy Littlewood method .
\par
Let
\begin{equation}
\label{eq:eight10}
V(x,Q) = \sum_{q\le Q}\sum_{a=1}^q \left|
\sum_{\substack{n\le x\\
n\equiv a\imod q}} r_3(n) - \Gamma(4/3)^3x\rho(q,a)q^{-3}
\right|^2
\end{equation}
Then the conclusion is
\begin{equation}
\label{eq:eight11}
V(x,Q)=Q\sum_{n\le x} r_3(n)^2 -A_1Qx+A_2Q^{\frac53}x^{\frac13} + E
\end{equation}
where
\[
E\ll x^{\frac{10}{9}+\varepsilon}\left(
\sum_{n\le x} r_3(n)^2
\right)^{\frac23} + Q^2(x/Q)^{\varepsilon}.
\]
Here
\[
A_1=\Gamma(4/3)^6 \sum_{q=1}^{\infty} q^{-6} \sum_{\substack{a=1\\
(a,q)=1}}^q |S_3(q,a)|^6
\]
as expected.  The surprising thing is that one has an asymptotic formula even though one cannot be certain of the size of the main term!

\Section{Thin sets}
\label{sec:nine}

\noindent So far for all the $a_n$ considered we have expected that
\[
\sum_{n\le x} a_n^2
\]
is roughly of order of magnitude $x$.  In Br\"udern \& Vaughan  \cite{BV21} we take instead $a_n=$
\[
r_2(n) = \card\{u,v\in\mathbb N^2:u^3+v^3=n\}.
\]
Now, by the usual lattice point arguments we have
\[
\sum_{n\le x} r_2(n) = Cx^{\frac23} + O(x^{\frac13}),\quad C=\frac{\Gamma(4/3)^2}{\Gamma(5/3)}
\]
and by a celebrated theorem of Hooley \cite{CH63}
\[
\sum_{n\le x} r_2(n)^2\sim 2Cx^{\frac23}.
\]
\par
It is natural to use the Hardy-Littlewood-method.  We consider
\[
V(x,Q) = \sum_{q\le Q}\sum_{a=1}^q \left|
\sum_{\substack{n\le x\\
n\equiv a\imod q}} r_2(n) - \frac{\rho(q,a)}{q^2} Cx^{2/3}
\right|^2
\]
and show that when $x^{3/5+\varepsilon}<Q\le x$
\[
V(x,Q) \sim 2CQx^{2/3}\sim Q\sum_{n\le x} r_2(n)^2.
\]
This is consistent with our overall philosophy since the major arcs are $\ll  Qx^{1/3}\log x$.
\par
It is noteworthy that the conclusion is deduced from a prior estimate for
\[
\sum_{q\le Q}\sum_{a=1}^q \left|
\sum_{\substack{n\le x\\
n\equiv a\imod q}} \left(
r(n) -\frac23 Cn^{-1/3}\mathfrak S(n;R)
\right)
\right|^2
\]
where $\mathfrak S(n;R)$, motivated by the idea of \S\ref{sec:seven}, is the truncated singular series
\[
\mathfrak S(n;R) = \sum_{r\le R} \sum_{\substack{b=1\\
(b,r)=1}}^r r^{-2} S_3(r,b)^2 e(-bn/r).
\]
\par
In Br\"udern \& Vaughan  \cite{BV22} we can also treat the case
\[
r(n)=\card\{u,v\in\mathbb N^2: u^k+v^l=n\}
\]
for various choices of $k< l$, namely
\[
k=2,\, l\ge 3\text{ and }k=3,l=4\text{ or }5.
\]
Let $\theta=\frac1k+\frac1l$,
\[
C=\frac{\Gamma(1+1/k)\Gamma(1+1/l)}{\Gamma(1+1/k+1/l)},
\]
$\rho(q,a)$ denote the number of solutions of the congruence $u^k+v^l\equiv a\imod q$ in $u$ and $v$, and let
\[
V(x,Q) = \sum_{q\le Q}\sum_{a=1}^q \left|
\sum_{\substack{n\le x\\
n\equiv a\imod q}} r(n) - \frac{\rho(q,a)}{q^2} Cx^{\theta}
\right|^2.
\]
Then $V(x,Q) \sim 2CQx^{\theta}$ for $x^{\theta-\eta}<Q\le x^{\theta}$ for some $\eta>0$ depending only on $k,l$.

\Section{Numbers with only small prime factors}
\label{sec:four}

\noindent A question which is significantly harder concerns the distribution of $y$-factorable numbers in arithmetic progressions.  For convenience we use the notation $\mathcal A(x,y)$, as in Vaughan  \cite{RV89}, to denote the set of positive integers not exceeding $x$ which have no prime factor exceeding $y$,
\begin{equation}
\label{eq:four1}
\mathcal A(x,y) = \{n\le x:p|n\Rightarrow p\le y\}.
\end{equation}
\par
This author deprecates the current prevalent use in this context of the barbarism ``smooth".  This has long established better uses in mathematics, and we will only use it here in those contexts.  The alternative term ``friable", whilst of reasonable provenance etymologically, in English really sounds like something which can be fried-up, like fish and chips.  It seems reasonable to describe numbers which can be factored into a product with none of the factors exceeding $y$ as being $y$-factorable.
\par
We have the familiar notation
\begin{equation}
\label{eq:four2}
\Psi_q(x,y) = \card\{n\in\mathcal A(x,y):(n,q)=1\},
\end{equation}
\begin{equation}
\label{eq:four2a}
\Psi(x,y)=\Psi_1(x,y),
\end{equation}
\begin{equation}
\label{eq:four3}
\Psi(x,y;q,a)=\card\{n\in\mathcal A(x,y):n\equiv a\imod q\}
\end{equation}
and
\begin{equation}
\label{eq:four4}
\Psi(x,y;\chi) = \sum_{n\in\mathcal A(x,y)} \chi(n)
\end{equation}
where $\chi$ is a Dirichlet character.  It is also useful to write
\begin{equation}
\label{eq:four5}
\mathcal A(y)=\mathcal A(\infty,y).
\end{equation}
\par
Now we are concerned with the behaviour of the variance
\begin{equation}
\label{eq:four6}
V(x,y,Q)\sum_{q\le Q} \sum_{a=1}^q  |\Psi(x,y;q,a) - \Upsilon(x,y;q,a)|^2
\end{equation}
where $\Upsilon$ is a suitable well behaved function.  When $(q,a)=1$ and $q$ is not too large in terms of $x$ it is known that
\[
\phi(q)^{-1}\Psi_q(x,y),
\]
is a good choice for $\Upsilon$ although this has the distinct disadvantage that it is discontinuous at $x$ and does not behave very smoothly at integral values of $\frac{\log x}{\log y}$.
\par
Comparisons with similar questions, as studied by Hooley in \cite{CH74}, ..., \cite{CH07}, by Br\"udern and Vaughan  \cite{BV21} and \cite{BV22}, by Br\"udern and Wooley \cite{BW11}, and by Vaughan  \cite{RV97}, ..., \cite{RV05}, would suggest that for a suitable $\Upsilon$ the above variance should be
\[
\approx Q\Psi(x,y).
\]
See also Harper \cite{AH12}.
\par
However, when $y=x$ it is clear that the best choice for $\Upsilon$ is $x/q$ and then
\[
V(x,x,Q)=\sum_{q\le Q} \sum_{a=1}^q \left|
\{-a/q\}-\{(x-a)/q\}
\right|^2 \ll Q^2,
\]
which is significantly smaller when $Q=o(x)$.  Moreover an  examination of the, apparently simple, case $\sqrt x<y<Q\le x$ already reveals some of the complications.  Thus the variability of $y$ clearly gives rise to additional features not normally occurring in generalizations of the Montgomery-Hooley theorem.  At least when $\sqrt x<y\le x$ we have
\[
\Psi(x,y;q,a) = \sum_{\substack{n\le x\\
n\equiv a\imod q}} 1 - \sum_{n\le x/y}\sum_{\substack{y<p\le x/n\\
np\equiv a\imod q}}1
\]
and so $V(x,y,Q)$ can be analysed by classical methods.  Thus generally we will concentrate on the case when $\log y=o(\log x)$.
\par
It transpires that we can imitate Hooley's method albeit in the variant utilised in Vaughan  \cite{RV03a} and are then able to leverage this to establish a corresponding bound for a more desirable main term, {\it viz}
\begin{equation}
\label{eq:four10x}
\Upsilon(x,y;q,a) = \frac1q\sum_{r|q}c_r(a)\sum_{s|r}\frac{\mu(r/s)}{\phi(r/s)}\Lambda_{r/s}^*(x/s,y)
\end{equation}
where $c_r(a)$ is the Ramanujan sum,
\begin{multline}
\label{eq:four15a}
\Lambda_q^*(x,y) = \frac{\phi(q)}{q}x\rho\left(
\frac{\log x}{\log y}
\right) \\
+ \int_0^{\infty} x\rho'\left(
\frac{\log x}{\log y}-v
\right)R_q(y^v) dv
\end{multline}
and $\rho$ is the Dickman function, defined on $\mathbb R$ to be continuous for $v\not=0$ and differentiable for $v\not=0,1$ and given by
\begin{equation}
\label{eq:four15b}
\rho(v)=\begin{cases}
0&(v<0),\\
1&(0\le v\le 1),
\end{cases}
\end{equation}
and
\begin{equation}
\label{eq:four15c}
v\rho'(v)+\rho(v-1)=0\quad(v>1).
\end{equation}
\par
The quality of the conclusion does not depend on whether or not there are exceptional Siegel zeros to the extent that the theory of the distribution of primes in arithmetic progressions does.  But the possibility of their existence does affect the organisation of the proof as well as the final conclusion.  However, all constants here, explicit or implicit, are effectively computable.
\par
Despite the circumlocutions the following theorem basically says that for a quite wide range of $x$ and $y$ we do indeed have an asymptotic formula for $V(x,y,Q)$ for a range of $Q$ that is normal in theorems of the Montgomery-Hooley kind.
\begin{theorem}[Vaughan \cite{RV23}]
\label{thm:four1}
There are small positive constants $c$ and $\kappa$ such that for each sufficiently small $\varepsilon>0$, when $x>x_0(\varepsilon)$ and $y$ satisfy
\begin{equation}
\label{eq:four10a}
\exp\left(
(\log\log x)^{\frac53+\varepsilon}
\right)\le y\le \exp\left(
\frac{\log x}{\log\log x}
\right)
\end{equation}
and $\xi$ and $K$ are defined by
\begin{equation}
\label{eq:four9}
\log\xi = \kappa\min\left(
\log y,\sqrt{\log x},\frac{\log x}{\log y}\log\log x
\right)
\end{equation}
and
\begin{equation}
\label{eq:four10}
K=\min\left(\xi^{c},\exp\big((\log y)^{\frac35-\varepsilon}\big)\right)
\end{equation}
the following holds.  Suppose that
\begin{equation}
\label{eq:four10b}
Q\le\Psi(x,y),
\end{equation}
and there is no $\xi$-exceptional zero.  Then
\begin{multline}
\label{eq:four10c}
\sum_{q\le Q}\sum_{a=1}^q\big|\Psi(x,y;q,a) - \Upsilon(x,y;q,a)\big|^2 \\
= Q\Psi(x,y)\bigg(1+O\Big(\exp\Big(-\frac{\log x}{2\log y}\log\frac{\log x}{\log y}\Big)\Big)\bigg) \\
+ O\big(\Psi(x,y)^2K^{-1}\big).
\end{multline}
If there is a $\xi$-exceptional zero, then the hypothesis has to be weakened to
\begin{equation}
\label{eq:four10d}
\exp\left(
(\log\log x)^{\frac53+\varepsilon}
\right)\le y\le \exp\left(
\frac{2\log x}{3(\log\log x)^4}
\right)
\end{equation}
and the error term has to be augmented by an additional term
\begin{multline}
\label{eq:four10e}
\Psi(x,y)^2(\log x)^5y^{-c/\log\log\log x} \\
+
\Psi(x,y)^2(\log x)^{15}q_{\xi}^{-1}x^{\beta_{\xi}-1}H\left(
{\textstyle\frac{\log x}{\log y}}
\right)^{-c}
\end{multline}
where $q_{\xi}$ is the conductor of the $\xi$-exceptional character and $H$ is given by \begin{equation}
\label{eq:four50}
H(u)=\exp\left(
\frac{u}{\big(\log(u+1)\big)^2}
\right).
\end{equation}
\end{theorem}

\Section{Conclusion}
\label{sec:ten}

\noindent It seems that it is possible to treat a wide range of questions of interest to number theorists.
All of the cases dealt with so far have
\[
\sum_{n\le x} a(n)^2 \gg x^{\lambda}
\]
where $\lambda>1/2$.  I have taken the second moment here so that one can include examples such as $a(n)=\mu(n)$.  By the way, I don't recall seeing it in the literature, but it is surely well known, and certainly easy to prove, that when $x(\log x)^{-A}\le Q=o(x)$ we have
\[
\sum_{q\le Q} \sum_{a=1}^q \left|
\sum_{\substack{n\le x\\
n\equiv a\imod q}} \mu(n)
\right|^2 \sim Q\sum_{n\le x}\mu(n)^2.
\]
Shparlinski has asked whether the methods described can treat the situation when $q$ is restricted to the primes, and the it is clear that the Hardy-Littlewood treatment can handle that.  Morever restrictions to squares, and more general $q$, are dealt with in Br\"udern and Wooley \cite{BW11}.

\Section{Questions}

\noindent {\bf Question 1.}  Is there a change of nature when
\[
\sum_{n\le x} a(n)^2 \ll x^{\lambda}
\]
with $\lambda<\frac12$?  Are there examples which behave differently?
Indeed, are there any examples!
\par
Hooley [1974] has conjectured that
\[
W(x,q) = \sum_{a=1}^q \Big|
\psi(x;q,a)-\frac{x}{\phi(q)}
\Big|^2 \sim x\log q
\]
and shown it p.p. $q$ with $\frac{Q}{2}<q\le Q$, $\frac{x}{(\log x)^{A}}<Q\le x$.
Fiorilli and Martin \cite{FM22}	have shown that this fails when $q$ is small with respect to $x$.  But surely Hooley would have been aware of this possibility and intended that $q$ is large, perhaps $q>x^{\delta}$.
\par
Both Hooley III [1975], and Friedlander and Goldston [1996] have extended the range for $Q$ on GRH. Also it can be shown Vaughan [2001] that if
\[
W(x,q) = \sum_{a=1}^q \Big|
\psi(x;q,a)-\frac{x}{\phi(q)}
\Big|^2
\]
and
\[
U(x,q) = x\log q - x\Big(
\gamma+\log 2\pi +\sum_{p|q}\frac{\log p}{p-1}
\Big),
\]
then there is an $F(y)$ a bit smaller than $y^{-1/2}$ so that.
\begin{align*}
M_k(x,Q) &= \sum_{Q/2<q\le Q} |W(x,q)-U(x,q)|^k \\
&\ll Qx^kF(x/Q)^k + Qx^k(\log x)^{-A}.
\end{align*}
\par
\noindent{\bf Question 2.}  What can be said about
\[
W(x,q) = \sum_{a=1}^q \left|
A(x;q,a)- \rho(q,a)\Psi(x)
\right|^2
\]
when we know a Montgomery-Hooley estimate for
\[
\sum_{q\le Q} \sum_{a=1}^q \left|
A(x;q,a)- \rho(q,a)\Psi(x)
\right|^2.
\]
Is there a much more general principle lurking here?  There are sequences $a_n$ known for which there is an asymptotic formula for $W(x,q)$ when, say, $x^{\theta}<q\le x$, see examples by Lau and Zhao \cite{LZ12}, Nunes \cite{RN15}.  Presumably this is a more general phenomenon and certainly there should be an estimate for almost all $q$.  I have not checked, but Hooley may have explored this to some extent in some of his many papers on Barban-Davenport-Halbestam.

\Section{Bibliography}
\label{sec:bib}

\end{document}